\documentclass[11pt]{article}
\usepackage{mathrsfs}

\input{amssymb.sty}
\usepackage{pb-diagram}
\begin{document}
\title{The Drazin inverse of the linear combinations of two
idempotents in the Banach algebra
\thanks{This project is supported by Natural Science
Found of China (10771191, 10471124 and 10771034).}}

\author{Shifang Zhang, \,\, Junde Wu \thanks{{Corresponding author:
wjd@zju.edu.cn}}}
\date{}
\maketitle

\begin{center}
\begin{minipage}{130mm}

\par{{\bf Abstract} In this paper, some Drazin inverse representations of the linear combinations of two
idempotents in Banach algebra are obtained.}

\vskip 0.2in


\vskip 0.2in


\end{minipage}
\end{center}

\vskip 0.3in

\section{\bf Introduction }
 Let $\mathscr{A}$ be a Banach$-*$ algebra with the unit $e$.
An element $P\in \mathscr{A}$ is said to be an idempotent if $P^2 =
P$ and a projection if  $P^2 = P=P^*$. The set
$\mathscr{P}(\mathscr{A})$ of all idempotents in $\mathscr{A}$ is
invariant under similarity, that is, if $P\in
\mathscr{P}(\mathscr{A})$ and $S\in \mathscr{A}$ is an invertible
element, then $S^{-1}PS$ is still an idempotent.

Let us recall that the Drazin inverse of $A \in \mathscr{A}$ is the
element $B\in \mathscr{A}$ (denoted by $A^D$) which satisfies
\begin{equation}BAB =B,\,\,\, AB =BA, \,\,\,A^ {k+1}B = A^ {k}\end{equation} for
some nonnegative integer $k$ ($\cite{1}$). The least such $k$ is the
index of $A$, denoted  by ind$(A)$. It is well-known that if $A$ is
Drazin invertible, then the Drazin inverse is unique and
$(aA)^D=\frac{1}{a}A^D$ for each nonzero scalar $a$. In particular,
for invertible operator $A$, the Drazin inverse $A^{D}$ coincide
with the usual inverse $A^{-1}$ and ind$(a) = 0$. The conditions (1)
are also equivalent to
\begin{equation} BAB =B,\,\,\,  AB =BA,\,\,\,  A-A^2B \makebox{\,\,is nilpotent.}\end{equation}
The Drazin inverse of an operator in $\mathscr{A}$ is similarly
invariant, that is, if $T$ is Drazin invertible and $S\in
\mathscr{A}$ is an invertible element, then $S^{-1}TS$ is still
Drazin invertible and $(S^{-1}TS)^D = S^{-1}T^DS.$ If $P\in
\mathscr{P}(\mathscr{A})$, it is easy to verify that $P^D=P$.

This paper is concerned with the Drazin inverses $(aP +bQ)^D$ of the
linear combinations of two idempotents in $\mathscr{A}$ for nonzero
scalars $a$ and $b$. In recent years, many authors paid much
attention to properties of linear combinations of  idempotents or
projections (see [2-7,9-14]). In $\cite{7}$, Deng has discussed  the
drazin inverses of the products and differences of two projections.
Motivated by this paper, A. B\"{o}ttcher and I. M. Spitkovsky wrote
$\cite{1}$ and in that paper they proved that  the Drazin
invertibility of the sum $P+Q$ of two projections $P$ and $Q$ is
equivalent to the Drazin invertibility of any linear combination $aP
+bQ$ where $ab\not=0, a+b\not=0.$ However, without some additional
conditions, it is difficult to discuss the Drazin invertibility  of
linear combinations of  two idempotents, even if the sum of them.
More recently, under some conditions,  Deng in $\cite{8}$ gave the
Drazin inverses of sums and differences of idempotents on the
Hilbert space. The methods used in $\cite{8}$ are the space
decompositions and operator matrix representations which are not
avail for general Banach$-*$ algebra, or general Banach algebra.

In this paper, by using the direct calculation methods, we obtained
some formulae for the Drazin inverse $(aP+bQ)^D$ of the linear
combinations of idempotents $P$ and $Q$ in Banach algebra
$\mathscr{A}$ under some conditions, we also study the index
ind$(aP+bQ)$.

\section{\bf Main results}

In this section, we always suppose that $\mathscr{A}$ is a Banach
algebra with the unit $I$, $aP+bQ$ is the linear combinations of two
idempotents $P$ and $Q$ in $\mathscr{A}$ with nonzero scalars $a$
and $b$. In order to prove $(aP+bQ)^D$ is Drazin invertible, it
follows from the definition of Drazin inverse that we only need to
find out some $M\in \mathscr{A}$ satisfies that

\begin{equation}  (aP+bQ)M =M(aP+bQ), M^2(aP+bQ) =M, (aP+bQ)^ {k+1}M = (aP+bQ)^
{k}\end{equation} for some nonnegative integer $k$.

\vskip0.1in

\vskip0.1in

\noindent {\bf Theorem 2.1.} Let $P$ and $Q$ be the idempotents in
Banach algebra $\mathscr{A}$ and $PQP=0$. Then $aP+bQ$ is Drazin
invertible for any nonzero scalars $a$ and $b$, ind$(a P+bQ)\leq 1$
and
$$(aP+bQ)^D=\frac{1}{a}P+\frac{1}{b}Q-(\frac{1}{a}+
\frac{1}{b})PQ-(\frac{1}{a}+\frac{1}{b})QP+(\frac{1}{a}+\frac{2}{b})QPQ.$$
Moreover, ind$(a P+bQ)=0$ if and only if $P+Q+QPQ=I+PQ+QP.$

\vskip0.1in

\noindent{\bf Proof.} We first prove that
$$(aP+Q)^D=\frac{1}{a}P+Q-(\frac{1}{a}+1)PQ-(\frac{1}{a}+1)QP+(\frac{1}{a}+2)QPQ.$$

\noindent For this, let
$M=\frac{1}{a}P+Q-(\frac{1}{a}+1)PQ-(\frac{1}{a}+1)QP+(\frac{1}{a}+2)QPQ.$
By the assumption that $PQP=0$, we have

$$\begin{array}{ll}&M(aP+Q)\\
=&(\frac{1}{a}P+Q-(\frac{1}{a}+1)PQ-(\frac{1}{a}+1)QP+(\frac{1}{a}+2)QPQ)(aP+Q)\\
=&[P+aQP-(a+1)P-(a+1)QP]+[\frac{1}{a}PQ+Q-(\frac{1}{a}+1)PQ-(\frac{1}{a}+1)QPQ+(\frac{1}{a}+2)QPQ]\\
=&P+Q-PQ-QP+QPQ\end{array}$$

\noindent and

$$\begin{array}{ll}
&(aP+Q)M\\
=&(aP+Q)(\frac{1}{a}P+Q-(\frac{1}{a}+1)PQ-(\frac{1}{a}+1)QP+(\frac{1}{a}+2)QPQ)\\
=&[P+aPQ-(a+1)PQ]+[\frac{1}{a}QP+Q-(\frac{1}{a}+1)QPQ-(\frac{1}{a}+1)QP+(\frac{1}{a}+2)QPQ]\\
=&P+Q-PQ-QP+QPQ.\end{array}$$

\noindent Therefore, $M(aP+Q)=(aP+Q)M.$ Moreover, a direct
calculation shows that

$$\begin{array}{ll}
&M(aP+Q)M\\
=&(P+Q-PQ-QP+QPQ)[\frac{1}{a}P+Q-(\frac{1}{a}+1)PQ-(\frac{1}{a}+1)QP+(\frac{1}{a}+2)QPQ]\\
=&[\frac{1}{a}P+PQ-(\frac{1}{a}+1)PQ]+[\frac{1}{a}QP+Q-(\frac{1}{a}+1)QPQ-(\frac{1}{a}+1)QP+(\frac{1}{a}+2)QPQ]\\
&-PQ-\frac{1}{a}QP-QPQ+(\frac{1}{a}+1)QPQ+QPQ
\\=&\frac{1}{a}P+Q-(\frac{1}{a}+1)PQ-(\frac{1}{a}+1)QP+(\frac{1}{a}+2)QPQ=M\end{array}$$
and $$\begin{array}{ll}
&M(aP+Q)^2\\
 =&\{\frac{1}{a}P+Q-(\frac{1}{a}+1)PQ-(\frac{1}{a}+1)QP+(\frac{1}{a}+2)QPQ\}
(aP+Q)^2\\
=&\{P+Q-PQ-QP+QPQ\}(aP+Q)\\
=&aP+aQP-aQP+PQ+Q-PQ-QPQ+QPQ=aP+Q.\end{array}$$ Thus, from (3) we
get that $(aP+Q)^D=M.$ So we have
$$\begin{array}{ll}
(aP+bQ)^D&=(b(\frac{a}{b}P+Q))^D\\
&=\frac{1}{b}(\frac{a}{b}P+Q)^D\\
&=\frac{1}{b}\{\frac{b}{a}P+Q-(\frac{b}{a}+1
)PQ-(\frac{b}{a}+1)QP+(\frac{b}{a}+2)QPQ\}
\\
&=\frac{1}{a}P+\frac{1}{b}Q-(\frac{1}{a}+
\frac{1}{b})PQ-(\frac{1}{a}+\frac{1}{b})QP+(\frac{1}{a}+\frac{2}{b})QPQ.
\end{array}$$
Moreover, since ind$(aP+Q)\leq1$ proved above and the fact that
ind$(aT)=$ind$(T)$ when $T$ is Drazin invertible, it follows that
ind$(aP+bQ)=$ind$(\frac{b}{a}P+Q)\leq1.$ In addition, a direct
calculation shows that $$(aP+bQ)^D(aP+bQ)=P+Q-PQ-QP+QPQ.$$ Note that
ind$(aP+bQ)=0$ if and only if $(aP+bQ)^D(aP+bQ)=I$, so
ind$(aP+bQ)=0$ if and only if $I=P+Q-PQ-QP+QPQ$. This completed the
proof.

\vskip0.1in

\noindent {\bf Theorem 2.2.} Let $P$ and $Q$ be the idempotents in
Banach algebra $\mathscr{A}$ and $PQP=P$. Then $aP+bQ$ is Drazin
invertible for any nonzero scalars $a$ and $b$, and

$$(aP+bQ)^D=\left\{\begin{array}{ll}\frac{a^2}{(a+b)^3}P+
\frac{1}{b}Q+\frac{ab}{(a+b)^3}(PQ+QP)+(\frac{b^2}{(a+b)^3}-\frac{1}{b})QPQ,&
\makebox{\,if\,\,} a+b\not=0;
\\\frac{1}{a}Q(P-I)Q,& \makebox{\,if\,\,}a+b=0.\end{array}\right.$$
Moreover, ind$(a P-aQ)\leq3$  and ind$(a P+bQ)\leq2$ when
$a+b\not=0$.

\vskip0.1in\noindent{\bf Proof.} Case (1). Suppose that $a+b\not=0.$
Firstly, we shall show that when $a\not=-1$, we have
$$(aP+Q)^D=\frac{a^2}{(a+1)^3}P+Q+\frac{a}{(a+1)^3}(PQ+QP)+(\frac{1}{(a+1)^3}-1)QPQ.$$
 To do this, let $M=\frac{a^2}{(a+1)^3}P+Q+\frac{a}{(a+1)^3}(PQ+QP)+(\frac{1}{(a+1)^3}-1)QPQ.$  By the
assumption that $PQP=P$, we have

 $$\begin{array}{ll}&M(aP+Q)\\
 =&\{\frac{a^2}{(a+1)^3}P+Q+\frac{a}{(a+1)^3}(PQ+QP)+(\frac{1}{(a+1)^3}-1)QPQ\}
(aP+Q)\\=&\frac{a^3}{(a+1)^3}P+aQP+\frac{a^2}{(a+1)^3}P+\frac{a^2}{(a+1)^3}QP+a(\frac{1}{(a+1)^3}-1)QP+
 \frac{a^2}{(a+1)^3}PQ\\
 &+Q+\frac{a}{(a+1)^3}PQ+\frac{a}{(a+1)^3}QPQ+(\frac{1}{(a+1)^3}-1)QPQ
\\
=&\frac{a^2}{(a+1)^2}P+Q+\frac{a}{(a+1)^2}(PQ+QP)+(\frac{1}{(a+1)^2}-1)QPQ\end{array}$$
and $$\begin{array}{ll}&(aP+Q)M\\
 =&(aP+Q) \{\frac{a^2}{(a+1)^3}P+Q+\frac{a}{(a+1)^3}(PQ+QP)+(\frac{1}{(a+1)^3}-1)QPQ\}
\\=&\frac{a^3}{(a+1)^3}P+aPQ+\frac{a^2}{(a+1)^3}PQ+\frac{a^2}{(a+1)^3}P+a(\frac{1}{(a+1)^3}-1)PQ+
 \frac{a^2}{(a+1)^3}QP\\
 &+Q+\frac{a}{(a+1)^3}QP+\frac{a}{(a+1)^3}QPQ+(\frac{1}{(a+1)^3}-1)QPQ
\\ =&\frac{a^2}{(a+1)^2}P+Q+\frac{a}{(a+1)^2}(PQ+QP)+(\frac{1}{(a+1)^2}-1)QPQ.\end{array}$$

Thus, \begin{equation} (aP+Q)M=M(aP+Q).\end{equation}

Since $$\begin{array}{ll}&M(aP+Q)^3\\
 =&\{\frac{a^2}{(a+1)^3}P+Q+\frac{a}{(a+1)^3}(PQ+QP)+(\frac{1}{(a+1)^3}-1)QPQ\}
(aP+Q)^3\\=&\{\frac{a^2}{(a+1)^2}P+Q+\frac{a}{(a+1)^2}(PQ+QP)+(\frac{1}{(a+1)^2}-1)QPQ\}
(aP+Q)^2\\
=&\{\frac{a^2}{(a+1)}P+Q+\frac{a}{(a+1)}(PQ+QP)+(\frac{1}{(a+1)}-1)QPQ\}
(aP+Q)\\ =&a^2P+Q+a(PQ+QP) \\ =&(aP+Q)^2, \end{array}$$ so,
\begin{equation} (aP+Q)^3M=(aP+Q)^2.\end{equation}
Moreover, by calculating, we get that   $$\begin{array}{ll}&M(aP+Q)M\\
 =&(\frac{a^2}{(a+1)^2}P+Q+\frac{a}{(a+1)^2}(PQ+QP)+(\frac{1}{(a+1)^2}-1)QPQ)\times
 \\& (\frac{a^2}{(a+1)^3}P+Q+\frac{a}{(a+1)^3}(PQ+QP)+(\frac{1}{(a+1)^3}-1)QPQ)
 \\=&\frac{a^4}{(a+1)^5}P+\frac{a^2}{(a+1)^3}QP+\frac{a^3}{(a+1)^5}QP+\frac{a^3}{(a+1)^5}P+(\frac{1}{(a+1)^2}-1)\frac{a^2}{(a+1)^3}QP+
 \\&\frac{a^2}{(a+1)^2}PQ+Q +\frac{a}{(a+1)^2}QPQ+\frac{a}{(a+1)^2}PQ+(\frac{1}{(a+1)^2}-1)QPQ+
 \\&\frac{a^3}{(a+1)^5}PQ+\frac{a}{(a+1)^3}QPQ+\frac{a^2}{(a+1)^5}QPQ+\frac{a^2}{(a+1)^5}PQ+(\frac{1}{(a+1)^2}-1)\frac{a}{(a+1)^3}QPQ+
 \\&\frac{a^3}{(a+1)^5}P+\frac{a}{(a+1)^3}QP+\frac{a^2}{(a+1)^5}QP+\frac{a^2}{(a+1)^5}P+(\frac{1}{(a+1)^2}-1)\frac{a}{(a+1)^3}QP+
  \\&\frac{a^2}{(a+1)^2}(\frac{1}{(a+1)^3}-1)PQ+\frac{a}{(a+1)^2}(\frac{1}{(a+1)^3}-1)QPQ+ (\frac{1}{(a+1)^3}-1)QPQ+
 \\& \frac{a}{(a+1)^2}(\frac{1}{(a+1)^3}-1)PQ+  (\frac{1}{(a+1)^2}-1)(\frac{1}{(a+1)^3}-1)QPQ

 \\=&(\frac{a^4}{(a+1)^5}+\frac{a^3}{(a+1)^5}+\frac{a^3}{(a+1)^5}+\frac{a^2}{(a+1)^5})P+Q+
 \\&(\frac{a^2}{(a+1)^3}+\frac{a^3}{(a+1)^5}+(\frac{1}{(a+1)^2}-1)\frac{a^2}{(a+1)^3}+\frac{a}{(a+1)^3}+\frac{a^2}{(a+1)^5}+(\frac{1}{(a+1)^2}-1)\frac{a}{(a+1)^3})QP+
 \\&(\frac{a^2}{(a+1)^2}+\frac{a}{(a+1)^2}+\frac{a^3}{(a+1)^5}+\frac{a^2}{(a+1)^5}+\frac{a^2}{(a+1)^2}(\frac{1}{(a+1)^3}-1)+\frac{a}{(a+1)^2}(\frac{1}{(a+1)^3}-1))PQ+
 \\&\{\frac{a}{(a+1)^2}+(\frac{1}{(a+1)^2}-1)+\frac{a}{(a+1)^3}+\frac{a^2}{(a+1)^5}+(\frac{1}{(a+1)^2}-1)\frac{a}{(a+1)^3}+\frac{a}{(a+1)^2}(\frac{1}{(a+1)^3}-1) +
\\&(\frac{1}{(a+1)^3}-1)+(\frac{1}{(a+1)^2}-1)(\frac{1}{(a+1)^3}-1)\}QPQ
\\=&(\frac{a^2}{(a+1)^3}P+Q+\frac{a^3+2a^2+a}{(a+1)^5}PQ+\frac{a^3+2a^2+a}{(a+1)^5}QP+\{\frac{a^2}{(a+1)^5}+\frac{1}{(a+1)^2}\frac{a}{(a+1)^3}+
\\&\frac{a}{(a+1)^2}\frac{1}{(a+1)^3}
+(\frac{1}{(a+1)^3}-1)+(\frac{1}{(a+1)^2}-1)\frac{1}{(a+1)^3}\}QPQ
 \\=&\frac{a^2}{(a+1)^3}P+Q+\frac{a}{(a+1)^3}(PQ+QP)+(\frac{1}{(a+1)^3}-1)QPQ=M.
\end{array}$$
That is, \begin{equation} M(aP+Q)M=M.\end{equation} It follows from
equations (4), (5) and (6) that  $aP+Q$ is Drazin invertible, $(a
P+Q)^D=M$ and ind$(a P+Q)\leq 2$ when $a\neq 1$.  Similar to the
disscussion in Theorem 2.1, when $a+b\not=0$, we have

$$\begin{array}{ll}
(aP+bQ)^D&=\frac{1}{b}(\frac{a}{b}P+Q)^D\\
&=\frac{1}{b}\{\frac{(\frac{a}{b})^2}{(\frac{a}{b}+1)^3}P+Q+\frac{\frac{a}{b}}{(\frac{a}{b}+1)^3}(PQ+QP)+(\frac{1}{(\frac{a}{b}+1)^3}-1)QPQ\}
\\&=\frac{a^2}{(a+b)^3}P+
\frac{1}{b}Q+\frac{ab}{(a+b)^3}(PQ+QP)+(\frac{b^2}{(a+b)^3}-\frac{1}{b})QPQ\end{array}$$
and ind$(a P+bQ)=$ind$(\frac{a}{b}P+Q)\leq 2$.

\vskip0.1in

Case (2). Suppose that $a+b=0.$ By calculating, we have
$$ (aP-aQ)\frac{1}{a}Q(P-I)Q= \frac{1}{a}Q(P-I)Q(aP-aQ)=Q-QPQ,$$

$$ (aP-aQ)(\frac{1}{a}Q(P-I)Q)^2= (Q-QPQ)\frac{1}{a}Q(P-I)Q=
\frac{1}{a}(QPQ-Q-QPQ+QPQ)=\frac{1}{a}Q(P-I)Q,$$ and

$$\begin{array}{ll}(aP-aQ)^4(\frac{1}{a}Q(P-I)Q)
&= (Q-QPQ)(aP-aQ)^3\\
&= a(QPQ-Q)(aP-aQ)^2\\
&= a^2(Q-QPQ)(aP-aQ)\\
 &=a^3(QPQ-Q)\\
 &=a^2(P-PQ-QP+Q)(aP-aQ)\\
  &=(aP-aQ)^3.\\  \end{array}$$
Therefore, $(aP-aQ)^D=\frac{1}{a}Q(P-I)Q$,
$(aP-aQ)^4(aP-aQ)^D=(aP-aQ)^3$ and ind$(a P-aQ)\leq 3.$ This
completed the proof.

\vskip0.1in

\noindent {\bf Remark 2.3.} Under the assumption of Theorem 2.2, we
have ind$(a P-aQ)= 3$ if and only if $P+QPQ\not=PQ+QP$. For this, we
only need to note that $(aP-aQ)^3(aP-aQ)^D=a^2(Q-QPQ)$ and
$(aP-aQ)^2=a^2(P-PQ-QP+Q).$

\vskip0.1in

\noindent {\bf Theorem 2.4.} Let $P$ and $Q$ be the idempotents in
Banach algebra $\mathscr{A}$ and $PQ=QP$. Then $aP+bQ$ is Drazin
invertible for any nonzero scalars $a$ and $b$, ind$(a P+bQ)\leq 1$
and
\begin{equation}(aP+bQ)^D=\left\{\begin{array}{ll}\frac{1}{a}P+\frac{1}{b}Q+(\frac{1}{a+b}-
\frac{1}{a}-\frac{1}{b})PQ,& \makebox{\,if\,\,} a+b\not=0;
\\\frac{1}{a}(P-Q),& \makebox{\,if\,\,}a+b=0.\end{array}\right.\end{equation}
Moreover, when $a+b\not=0$, $ind(a P+bQ)=0$ if and only if
$P+Q=I+PQ$; while ind$(a P-aQ)=0$ if and only if $P+Q=I+2PQ$.

\vskip0.1in

\noindent{\bf Proof.} We first prove that when $a\not=-1,$
$$(aP+Q)^D=\frac{1}{a}P+Q+(\frac{1}{a+1}-\frac{1}{a}-1)PQ.$$
For this, let $M=\frac{1}{a}P+Q+(\frac{1}{a+1}-\frac{1}{a}-1)PQ$. By
the assumption that $PQ=QP$, a direct calculation shows that
$$\begin{array}{ll}(aP+Q)M&=M(aP+Q)\\
 &=(\frac{1}{a}P+Q+(\frac{1}{a+1}-\frac{1}{a}-1)PQ)(aP+Q)\\
 &=P+aPQ+(\frac{a}{a+1}-1-a)PQ+\frac{1}{a}PQ+Q+(\frac{1}{a+1}-\frac{1}{a}-1)PQ\\
  &=P+Q-PQ.\\  \end{array}$$ Moreover, it is easy to check  that $$\begin{array}{ll}(aP+Q)^2M&=(aP+Q)(P+Q-PQ)\\
  &=aP+aPQ-aPQ+PQ+Q-PQ
\\&=aP+Q\end{array}$$
and
$$\begin{array}{ll}&M(aP+Q)M\\
&=(P+Q-PQ)(\frac{1}{a}P+Q+(\frac{1}{a+1}-\frac{1}{a}-1)PQ)\\
&=\frac{1}{a}P+\frac{1}{a}PQ-\frac{1}{a}PG+PQ+Q-PQ+(1+1-1)(\frac{1}{a+1}-\frac{1}{a}-1)PQ
\\&=\frac{1}{a}P+Q+(\frac{1}{a+1}-\frac{1}{a}-1)PQ=M.\end{array}$$
So, $$(aP+Q)^D=\frac{1}{a}P+Q+(\frac{1}{a+1}-\frac{1}{a}-1)PQ.$$ If
$a+b\not=0$, then
$$\begin{array}{ll}
(aP+bQ)^D&=(b(\frac{a}{b}P+Q))^D\\
&=\frac{1}{b}(\frac{a}{b}P+Q)^D\\
&=\frac{1}{b}\{\frac{b}{a}P+Q+(\frac{b}{a+b}-\frac{b}{a}-1)PQ\}
\\&=\frac{1}{a}P+\frac{1}{b}Q+(\frac{1}{a+b}-
\frac{1}{a}-\frac{1}{b})PQ.\end{array}$$ Moreover, we can show that
ind$(aP+bQ)\leq1$ and when $a+b\not=0$,
$$(aP+bQ)^D(aP+bQ)=P+Q-PQ.$$
So, ind$(aP+bQ)=0$ if and only if $ I=P+Q-PQ$.

\vskip0.1in

On the other hand, note that $PQ=QP$, so we have
$$(P-Q)^2= P+Q-2PQ \makebox{\,\, and \,\,} (P-Q)^3= P-Q,$$
this implied that $(P-Q)^D= P-Q$. Thus, when $a+b=0$, we have
$(aP+bQ)^D=\frac{1}{a}(P-Q)$ and ind$((aP+bQ)^D)\leq1$. It is clear
that ind$(aP-aQ)=0$ if and only if $P+Q=I+2PQ$. This completed the
proof.

\vskip0.1in

Noting that  $PQP=Q$ implies that $Q=QP=PQ$, so, it follows from
Theorem 2.4 immediately:

\vskip0.1in\noindent {\bf Corollary 2.5.} Let $P$ and $Q$ be the
idempotents in Banach algebra $\mathscr{A}$ and $PQP=Q$. Then
$aP+bQ$ is Drazin invertible for any nonzero scalars $a$ and $b$,
ind$(a P+bQ)\leq 1$ and
$$(aP+bQ)^D=\left\{\begin{array}{ll}\frac{1}{a}P+(\frac{1}{a+b}-
\frac{1}{a})Q,& \makebox{\,if\,\,} a+b\not=0;
\\\frac{1}{a}(P-Q),& \makebox{\,if\,\,}a+b=0.\end{array}\right.$$

\vskip0.1in

\noindent {\bf Remark 2.6.} (1). It follows from Corollary 2.5 that
if $PQP=Q$, then $(P-Q)^D=P-Q$. Moreover, we can prove that
 $(P-Q)^D=P-Q$ if and only if  $PQP=QPQ$.

(2). Our results recovered most of the main conclusions in
$\cite{8}$, but our methods are very different from the methods used
in $\cite{8}$, in particular, the methods used in $\cite{8}$ cannot
obtain any information about the Drazin index.

\vskip 0.9in
 The group inverse of $A \in \mathscr{A}$  ([16-19]) is the
element $B\in \mathscr{A}$ (denoted by $A^g$) which satisfies
\begin{equation}BAB =B,\,\,\, AB =BA, \,\,\,ABA = A.\end{equation}
 Obviously,  $A$ has group inverse if and only if
 $A$ has Drazin inverse with  ind$(A)\leq1$.

\vskip 0.2in Before  giving the  revised versions of theorems 3.2
and 3.3 in $\cite{15}$, let us see the  following two interesting
counter-examples.

\noindent{\bf Example 2.7} Let
$A=\left(\begin{array}{cc}S&0\\0&0\end{array} \right)\in B(l_2\oplus
l_2)$ and
 $B=\left(\begin{array}{cc}0&0\\T&0\end{array} \right)\in
B(l_2\oplus l_2)$ with $S$ and $T$ in $B(l_2)$  such that $TS\neq0$.
Considering operator

$$P=\left(\begin{array}{cc}I&0\\0&0\end{array} \right)\in B(H_2\oplus
H_2),\,\,\,\,\,\,Q=\left(\begin{array}{cc}I&A\\B&0\end{array}
\right)\in B(H_2\oplus H_2),$$ where $H_2=l_2\oplus l_2$.

Direct calculations shows that $$BA\neq0,\,\,\,(BA)^2=AB=0.$$

\noindent Hence we have $P^2=P, Q^2=Q, PQP=P$. From Theorem 2.2 , we
know that $P+Q$ has Drazin inverse and $(P+Q)^D=\frac{1}{8}P+
Q+\frac{1}{8}(PQ+QP)-(\frac{7}{8})QPQ$. Hence
$(P+Q)-(P+Q)^2(P+Q)^D=(P+Q)-(\frac{1}{2}P+Q+\frac{1}{2}(PQ+QP)-\frac{1}{2}
 QPQ)=\frac{1}{2}(P+QPQ)-\frac{1}{2}(PQ+QP)=\frac{1}{2}\left(\begin{array}{cc}0&0\\0&BA\end{array}
\right)\neq0$, which implies that  ind$(P+Q)>1$. Together this with
Theorem 2.2, it is clear that ind$(P+Q)=2.$ So
 the group inverse $(P+Q)^g$ of
$P+Q$ does not exist.

\noindent{\bf Example 2.8}   Define operators $p$ and $q$ in
$B(\mathbb{C}^5)$ by
$p=\left(\begin{array}{ccccc}1&0&0&0&0\\0&1&0&0&0\\0&0&0&0&0\\0&0&0&0&0\\0&0&0&0&0\end{array}
\right)$ and
$q=\left(\begin{array}{ccccc}1&0&0&0&0\\0&1&0&0&0\\0&0&0&0&0\\0&1&0&0&0\\0&0&0&0&1\end{array}
\right)$, respectively. Obviously,
$$p^2=p, q^2=q, pqp=p=pq.$$
This means that $p$ and $q$ are idempotents in $B(\mathbb{C}^5)$.
Then it results from Theorem 2.2 that $(p-q)^D=q(p-1)q$. But a
direct calculation shows that
$(p-q)^2(p-q)^D=qpq-q=\left(\begin{array}{ccccc}0&0&0&0&0\\0&0&0&0&0\\0&0&0&0&0\\0&0&0&0&0\\0&0&0&0&-1\end{array}
\right)\neq p-q$, this mean that  ind$(p-q)> 1$. So
 the group inverse $(p-q)^g$ of
$p-q$ does not exist.

\vskip 0.1in

 The above two examples  illustrate not only Theorem 3.2, but also part (ii) of Theorem 3.3 in $\cite{15}$
 are not always true. Now we present the  modified versions as follows

\vskip 0.2in

 \noindent{\bf Theorem 3.2$'$} Let $P$ and $Q$ be the idempotents in
Banach algebra $\mathscr{A}$ and $PQP=P$  Then

(i)$(P+Q)^D=\frac{1}{8}P+Q+\frac{1}{8}(PQ+QP)-(\frac{7}{8})QPQ$,

 (ii) $(P-Q)^D=Q(P-1)Q$,

 (iii)  $P+Q$ has group inverse if and only if $
P+QPQ=PQ+QP$ ,

 (iv) $P-Q$ has group inverse if and only if $ P=QPQ$.

 \noindent{\bf Proof.}  Since  the results of part (i) and  part
 (ii) is a special case of Theorem 2.2 , it suffice  to show part (iii) and  part
 (iv). For this, we only need to note that
$(P+Q)-(P+Q)^2(P+Q)^D=\frac{1}{2}(P+QPQ-PQ-QP)$ and that
$(P-Q)-(P-Q)^2(P-Q)^D=P-QPQ$, which can be obtained by direct
calculations.  This completed the proof.

\vskip0.1in

\noindent {\bf Theorem 3.3$'$} Let $P$ and $Q$ be the idempotents in
Banach algebra $\mathscr{A}$ and $PQP=PQ$. Then
$$(P+Q)^g=P+Q-2QP-\frac{3}{4}PQ+\frac{5}{4}QPQ,$$
$$(P-Q)^D=P-Q-PQ+QPQ.$$

Moreover, ind$(P-Q)\leq2$ and  $P-Q$ has group inverse if and only
if $PQ=QPQ$.

\vskip0.1in\noindent{\bf Proof.} Since the group inverse of $P +Q$
can by checked directly, its proof is omitted. Now let
$M=P-Q-PQ+QPQ$. By direct calculations we have that \begin{equation}
M(P-Q)M=M, (P-Q)^2M=M,\end{equation} and that
$$ (P-Q)^3M=(P-Q)^2=(P-Q)M=M(P-Q)=P-PQ-QP+Q.$$

\noindent This implies that $(P-Q)^D=P-Q-PQ+QPQ$ and that
ind$(P-Q)\leq2$. In this case, from equation (9) and the definition
of group inverse, we know that $P-Q$ has group inverse if and only
if $(P-Q)^2(P-Q)^D=(P-Q)=(P-Q)^D=P-Q-PQ+QPQ.$ This completed the
proof.

\end{document}